\documentclass[12pt,draftcls,onecolumn]{IEEEtran}
\usepackage{mathrsfs,amsmath,latexsym,amssymb}

\newtheorem{thm}{Theorem}[section]

\newcommand{\argmax}{\mathop{\mbox{\rm arg\,max}}}
\newcommand{\argmin}{\mathop{\mbox{\rm arg\,min}}}
\newcommand{\real}{\mathbb R} 
\begin{document}
\sloppy

\title{Top Feasible-Arm Subset Identification in Constrained Multi-Armed Bandit with Limited Budget}
\author{Hyeong Soo Chang
\thanks{H.S. Chang is with the Department of Computer Science and Engineering at Sogang University, Seoul 121-742, Korea. (e-mail:hschang@sogang.ac.kr).}%
}

\maketitle
\begin{abstract}
We present an algorithm, ``constrained successive accept or reject (CSAR),"
for the problem of identifying the subset of top feasible-arms from a given 
finite set of arms with the limited sampling-budget equal to a given time-horizon 
when the sequential dynamics 
of the arms follows the model of a constrained multi-armed bandit.
We provide a finite-time upper bound on the probability of the incorrect identification 
by CSAR that converges to zero with an exponential rate in the sampling-budget.
\end{abstract}


\section{Introduction}

We consider the problem of identifying the subset of top feasible-arms from a given finite set of arms with the limited sampling-budget equal to a given time-horizon when
the sequential dynamics of the arms follows the model of a constrained multi-armed bandit (CMAB). 
The model is a special case of constrained Markov decision process (see, e.g.,~\cite{changcmab}~\cite{changindex}) as classical stochastic MAB is known to be a special case of Markov decision process.

Our goal is to present an algorithm, ``constrained successive accept or reject (CSAR)," extending ``successive accept or reject (SAR)" by Bubeck \emph{et al.}~\cite{bubeck} designed for the unconstrained problem
providing a finite-time upper bound, given with a problem character constant, on the probability of the incorrect identification by
CSAR that converges to zero with an exponential rate in the sampling-budget.

At each discrete time, a single arm is selected and played by CSAR. Once played, not only a sample reward but also a sample cost
are independently drawn from the unknown reward distribution and the unknown cost distribution associated with
the arm, respectively. By a given time-horizon, a set of arms is returned as a solution by CSAR, which is possibly empty.

The values of the mean and the variance of each distribution for reward and cost associated with each arm are unknown. However, we assume that the reward means are all different among the arms and the cost means are also all different among the arms.

It should be noted that CSAR maintains salient structure of SAR and also follows the main idea of the performance analysis \emph{but with subtle and important modifications} in order to handle the constrained setting while addressing some issues that cannot be ignored about the algorithmic incompleteness of SAR and the performance proof (see, subsection~\ref{sar}).

There exists a great body of literature on problems about finding top (feasible) designs from a finite set of designs in the order of correct ranking (or their variants) in the area of ranking and selection under (constrained) simulation optimization contexts (see, e.g.,~\cite{li}~\cite{xiao}~\cite{zhang}~\cite{pasupathy} and the reviews and the references therein).
But most notably the problem setting is different from the bandit setting:
once a particular design is picked, \emph{multiple} reward (and/or cost) samples can be drawn at a single time-step from a distribution or a reward (batch) sample (that possibly depends on the states of all alternative designs) is obtained \emph{after} the total simulation-replications are \emph{used up}.
Furthermore, those samples in general are assumed to be \emph{normally} distributed. 
Under the normality (and possibly pre-knowledge on some problem parameters) then, the main interest is in an optimal sampling allocation of assigning the number of simulation replications to each design for a given performance measure.
Some iterative heuristic algorithms that incrementally update the allocations over the designs have been developed to approximate the optimal sampling-allocation (see, e.g.,~\cite{li}~\cite{xiao}~\cite{zhang} and the references therein).

Related constrained problems in MAB settings under various contexts have been also studied
(see, e.g.,~\cite{katz}~\cite{hou}~\cite{sinha} and the references therein).
However, the dynamics of the arms in those works is different from that of CMAB. Each arm's sample used for measuring its performance and its feasibility is generated from a \emph{single} distribution.
For the problem of identifying a \emph{single} best feasible-arm in CMAB, Chang~\cite{changcmab} provides a convergent stochastic algorithm but it achieves only \emph{asymptotic} optimality.
An index-based deterministic algorithm has been developed recently where it achieves an exponential convergence rate in the horizon size~\cite{changindex} but computing the indices requires knowing the mean of the best feasible-arm \emph{in advance}.

\section{Top Feasible-Arm Subset Identification: Algorithm}

\subsection{Preliminaries}

A finite set $A$ of arms is given with $|A|>1$.
At each discrete time $t\geq 1$, a single arm $a_t$ in the \emph{active}-arm set $A^t \in 2^A\setminus \emptyset$ is selected and played (or pulled).
Then a sample reward $X_{a_t,t}$ in $[0,1]$ and a sample cost $Y_{a_t,t}$ in $[0,1]$ are independently drawn from the unknown reward distribution and the unknown cost distribution associated with $a_t$, respectively. 
(For simplicity, we assume that the distributions are all supported on $[0,1]$.)
For any $a\in A$ and $t\neq t'$, $E[X_{a,t}] = E[X_{a,t'}]$ and $E[Y_{a,t}] = E[Y_{a,t'}]$.

Let $C:A\rightarrow \real$ be given such that $C(a)= E[Y_{a,t}], a\in A,$ and $\mu:A \rightarrow \real$ be given such that $\mu(a)=E[X_{a,t}], a\in A,$ with any fixed $t$.
Given a constant $\tau > 0$, $a$ in $A$ is \emph{feasible} if $C(a) \leq \tau$ and let
\[A_f=\{a| C(a) \leq \tau, a\in A \}.
\] It should be noted that the size of $A_f$ is \emph{unknown} except that $0\leq |A_f| \leq |A|$.
Let $\{\mu_{(i)}| i = 1,...,|A_f|\}$ be an \emph{ordered} set obtained by permutation of $\{\mu(a)| a \in A_f\}$ in the decreasing order of $\mu$-values when $A_f \neq \emptyset$, where the subscript $(i)$ in $\mu_{(i)}$ refers to the $i$th \emph{largest} element in the set.

The problem is, if $A_f = \emptyset$, to find the empty set as the solution, else if $A_f \neq \emptyset$ and $1\leq m < |A_f|$, to find a subset $A_f(m) =\{\mu_{(i)} | i=1,...,m\}$ of $A_f$ that contains the top-$m$ feasible-arms,
or else if $A_f \neq \emptyset$ and $|A_f| \leq m \leq |A|$, to find $A_f$ as the solution, by playing the bandit up to a fixed finite time-horizon $H$.
Because $H$ is fixed, the number of the plays or the sampling-budget is \emph{limited} by $H$.
(For the problem of selecting the top-$m$ arms in \emph{unconstrained} MABs, because we know $|A|$, $1 \leq m \leq |A|-1$ in general. However, in our case, not only we do not know $|A_f|$ but also it is possible that $|A_f|=|A|$. Therefore, we consider that $1\leq m\leq |A|$.)

Given a function $g:S \rightarrow \real$ for a given nonempty set $S$, we define a bijection $\sigma_{g}[S]:S\rightarrow \{1,...,|S|\}$ such
that for $s\in S$, $\sigma_{g}[S](s) = i$ if and only if $g(s)$ is the $i$th largest element in $S$.
In other words, $\sigma_{g}[S](s)$ is the \emph{rank} of $s$ in terms of $g$-value among the elements of $S$.
Therefore, $g(\sigma_g[S]^{-1}(1)) = \max_{s\in S} \{g(s)\}$ and $\sigma_g[S]^{-1}(1) = \argmax_{s\in S} \{g(s)\}$. 
We assume that if $\max_{s\in S} \{g(s)\}$ is achieved by more than one elements in $S$, then $\argmax_{s\in S} \{g(s)\}$ breaks the tie by the arm whose rank in $g$ is the highest, i.e.,
\[\argmax_{s\in S} \{g(s)\} = \argmin_{s'\in S}\{ \sigma_{g}[S](s')| g(s') = \max_{s\in S} \{g(s)\}, s' \in S \}.\]

Furthermore, given a function $g$ defined with a nonempty set $S$ as its domain, for any nonempty $S'\subseteq S$,
(with abusing the notation) we define the \emph{ranking restricted with $S'$} (of $S$) as a bijection $\sigma_{g}[S']:S'\rightarrow \{1,...,|S'|\}$ such that for $s\in S'$, $\sigma_{g}[S'](s)$ is the rank of $s$ in terms of $g$-value \emph{among the elements in} $S'$ (instead of $S$).

Define $\Delta_c: A \rightarrow \real$ such that for $a\in A$, $\Delta_c(a) = |C(a) - \tau|$.
In the sequel, we assume that $C(a)\neq \tau$ for all $a\in A$. That is, $\Delta_c(a)>0$ for all $a\in A$.

For a nonempty $A_f$ and $1\leq m < |A_f|$, define $\Delta: A \rightarrow \real$ such that for $a\in A_f$,
\[ \Delta(a) = \begin{cases}
    \mu_{(\sigma_{\mu}[A_f](a))} - \mu_{(m+1)} & \text{if } \sigma_{\mu}[A_f](a) \leq m \\
    \mu_{(m)} - \mu_{(\sigma_{\mu}[A_f](a))} & \text{otherwise}.
\end{cases}
\]
and for all $a\in A\setminus A_f$, $\Delta(a) = \min_{a\in A_f} \Delta(a)$.
The function $\Delta$ is \emph{not} defined when $A_f = \emptyset$ or $m \geq |A_f|$.
Note that from the definition of $\Delta$, either both $\sigma_{\Delta}[A_f]^{-1}(1)$ and $\sigma_{\Delta}[A_f]^{-1}(|A_f|)$ or 
one of them \emph{only} achieve(s) $\max_{a \in A_f} \Delta_{(a)}$.

We proceed an empirical-estimate of $\Delta$ obtained by CSAR.
CSAR works with the \emph{phase} $k=1,...,|A|$ and the parameter $N_k$ for each phase $k$.
Let $t^k = |A|N_1 + (|A|-1) N_2 + \cdots + (|A|-k+1) N_k$ for $k=1,...,|A|$.
The phase $k$ corresponds to the time steps in $\{t^{k-1}+1,t^{k-1}+2,t^{k-1}+3,...,t^k\}$, where $t^0 = 0$.
It is required that the values of $N_k$'s for $k=1,...,|A|$ should be chosen to satisfy $t^{|A|}\leq H$.

At each phase $k$,
CSAR is associated with a set $A^k \subseteq A$ that contains $|A|-k+1$ active arms obtained by deactivating
exactly one arm from $A^{k-1}$ with $A^1=A$. 
$N_k$-samples are
drawn for each arm in $A^k$ ($N_k \times |A^k|$ in total) over the time steps in the phase $k$.
In other words, 
$(|A|-k+1)N_k$-samples are 
drawn over the time steps in the phase $k$ in an arbitrary order of the arms in $A^k$ such that
\[\bigcup_{n=1,...,k}\bigcup_{a\in A^n} N_a(t^n) = \{1,...,t^k\},
\] where $N_a(t) = \{n|a_n =a, n=1,...,t\}$ for $a\in A$ referring to the set of the time steps when the arm $a$ was pulled up to $t$.

For a fixed $k \in \{1,...,|A|\}$ and a given nonempty set $A^k\subseteq A$, define
$\hat{C}^k:A^k \rightarrow \real$ that provides the sample mean of the cost samples up to time $t^k$ such that
for $a\in A^k$,
$\hat{C}^k(a) =  \frac{1}{|N_a(t^k)|} \sum_{n \in N_a(t^k)} Y_{a,n}$ and define
$\hat{\mu}^k: A^k \rightarrow \real$ that provides the sample mean of the reward samples up to time $t^k$ such that for $a\in A^k$,
$\hat{\mu}^k(a) =  \frac{1}{|N_a(t^k)|} \sum_{n \in N_a(t^k)} X_{a,n}$.

CSAR finds the set $A_f^k$ of empirically-feasible arms in $A^k$ from $\hat{C}^k$-values 
and works with a counting variable $m_k$ used for tracing the number of the arms that have been selected as top-$m$ arms at the phase $k$. Only when $A_f^k \neq \emptyset$ and $m_k < |A_f^k|$, 
an empirical-estimate of $\Delta$ is used for determining acceptance or rejection:
For a fixed $k \in \{1,...,|A|\}$, and a given nonempty set $A_f^k\subseteq A^k$, 
and a given $m_k \in \{1,...,|A_f^k|-1\}$, we define
$\hat{\Delta}^k: A_f^k \rightarrow \real$ such that for $a\in A_f^k$
\begin{eqnarray}
\label{deltahat}
\hat{\Delta}^k(a) = 
\begin{cases}
\hat{\mu}^k_{(\sigma_{\hat{\mu}^k}[A_f^k](a))} - \hat{\mu}^k_{(m_k+1)} & \mbox{if } \sigma_{\hat{\mu}^k}[A_f^k](a) \leq m_k, \\
\hat{\mu}^k_{(m_k)} - \hat{\mu}^k_{(\sigma_{\hat{\mu}^k}[A_f^k](a))} & \mbox{otherwise},
\end{cases} 
\end{eqnarray} where
$\{\hat{\mu}^k_{(i)}| i = 1,...,|A_f^k|\}$ is an ordered set obtained by permutation of $\{\hat{\mu}^k(a)| a \in A_f^k\}$ in the non-increasing order of $\hat{\mu}^k$-values.

The function $\hat{\Delta}^k$ is \emph{not} defined (i.e., not used in CSAR) when $A_f^k = \emptyset$ or $m_k \geq |A_f^k|$.
Note that similar to $\Delta$, either both $\sigma_{\hat{\Delta}^k}[A_f^k]^{-1}(1)$ and $\sigma_{\hat{\Delta}^k}[A_f^k]^{-1}(|A_f^k|)$ or one of them achieve(s) $\max_{a \in A_f^k} \hat{\Delta}^k_{(a)}$.

In the next subsection, we describe SAR by using the terms introduced here.

\subsection{Successive Accept or Reject: Issues}
\label{sar}

Assume that $A_f=A$ in this subsection to discuss the \emph{unconstrained} problem of identifying the subset of top-$m$ arms in $A$ with $1\leq m \leq |A|-1$. The pseudocode of SAR is given below.
SAR begins the loop with the set of active arms $A^k$, and selects an arm $d_k$ to deactivate from $A^k$ (from $\hat{\Delta}^k$ defined with $A_f^k=A^k$ in~(\ref{deltahat})). The counting variable $m_k$ takes the value of how many arms need to be accepted more as top-$m$ arms:
\\
\\
\noindent\textbf{Successive Accept or Reject (SAR)}
\begin{itemize}
\item[1.] \textbf{Initialization:} Set $T_i=$Nil for $i=1,...,m$, $A^1 = A$, $m_1=m$, and $k=1$.\\
Select $N_i$ for $i=1,...,|A|-1$ such that $\sum_{k=1}^{|A|-1} N_k (|A|-k+1) \leq H$.
\item[2.] \textbf{Loop while $k\leq |A|-1$}
\begin{itemize}
\item[2.1] Play each $a \in A^k$ $N_k$ times.
\item[2.2] \textbf{Deactivation:} Choose any arm in $\{a\in A^k| \hat{\Delta}^k(a) = \max_{s\in A^k}  \hat{\Delta}^k(s) \}$ as $d_k$ and $A^{k+1} = A^k \setminus \{d_k\}$.
\item[2.3] \textbf{Accept or Reject:} If $d_k = \sigma_{\hat{\mu}^k}[A^k]^{-1}(1)$, then $T_{m-m_k+1} = d_k$ and $m_{k+1} = m_k - 1$. Set $k \leftarrow k+1$.
\end{itemize}
\item[3.] \textbf{Output:} Return $\{T_1,...,T_m\}$.
\end{itemize}
\vspace{0.5cm}

As we can see, SAR is simple. However, the algorithmic description of SAR is still incomplete in that the tie-condition is \emph{not} handled properly.
When SAR selects an arm to deactivate in the step 2.2, it \emph{breaks ties arbitrarily} (as explicitly stated in the pseudocode of SAR in~\cite{bubeck}).

Because $d_k$ is set to be \emph{any} element in $\{ a | \Delta^k(a) = \max_{s\in A^k}  \Delta^k(s) \}$, SAR can \emph{reject incorrectly}.
Suppose that $k=|A|-1$, $|A^k|=2$, and $m_{|A|-1}=1$ and the highest empirical-gap $\Delta^k$ is achieved by \emph{both} $\sigma_{\hat{\mu}^k}[A^k]^{-1}(1)$ and $\sigma_{\hat{\Delta}^k}[A^k]^{-1}(|A^k|)$. 
If we break the tie by
$\sigma_{\hat{\Delta}^k}[A^k]^{-1}(|A^k|)$, then because the arm does not satisfy the acceptance condition, it
rejects $\sigma_{\hat{\Delta}^k}[A^k]^{-1}(|A^k|)$. SAR finishes \emph{without setting} $T_m$ making the output $\{T_1,...,T_m\}$ fail.
SAR should have accepted $\sigma_{\hat{\mu}^k}[A_f^k]^{-1}(1)$ to correctly finish.
 
Because this tie-situation can happen at any phase $k$, SAR might not correctly reject at any $k$.
Indeed, Bubeck \emph{et al.}~missed incorporating their tie-breaking rule into their correctness proof. (For example, the statement ``$a_j = \sigma(1)$" made in the proof of Theorem 1 in~\cite{bubeck} should match their tie-breaking rule.)
A fix is simple. SAR can resolve this issue by the consistent selection such that the tie is broken with the highest rank in the argument set of the maximum. In fact, this has been already specified in the previous subsection for the argmax operator and CSAR uses this rule.
Even if the remedy is minor, the issue cannot be ignored in order for the proof~\cite{bubeck} to be correct.

Next, SAR does not have any \emph{exit condition}. SAR stops only when $k > |A|-1$.
Even if $m_{k+1}$ is set to zero at some $k$, SAR still runs at $k+1$ if $k+1 \leq |A|-1$.
This can cause $T_{m+1}$ to be set with $d_{k+1}$ and $m_{k+2} = - 1$, and cause $T_{m+2}$ to be set, etc. When $m_k=0$ at the beginning of the loop, SAR should exit the loop.
Furthermore, if $|A| \gg m$, it is possible that SAR runs ``on empty" for many iterations.

Bubeck \emph{et al.}'s idea of the proof about bounding the probability of incorrect identification can be followed.
But their whole proof is based on an \emph{erroneous definition} of the conditioning event
(notated by $\xi$ in~\cite{bubeck}), which
is actually \emph{impossible}. This necessarily affects the validity of the proof.
At each phase of SAR, \emph{only active} arm's samples are drawn. Only the sample means 
of the active arms at the phase $k$, but not all arms in $A$, are updated with the samples drawn up to
the phase $k$.
Because the probability of $\bar{\xi}$ is just one, the upper bound given in~\cite{bubeck} on the probability of error by SAR cannot be implied.

\subsection{Constrained Successive Accept or Reject}

At each phase $k=1,...,|A|$, CSAR begins the loop with the set of active arms $A^k$ as in SAR, from which the set of active empirically-feasible arms $A^k_f$ is obtained.
The variable $m_k$ has the value of how many arms need to be accepted more as top-$m$ feasible arms.

After each arm in $A^k$ is played $N_k$ times (in an arbitrary order), CSAR 
finds $A_f^k$ from $\hat{C}^k$-values,
which contains the arms in $A^k$ that empirically satisfies the feasibility condition.
Then it selects $d_k$ to be deactivated in $A_f^k$ in three ways: If $|A^k_f|>m_k$, then
the highest rank arm that achieves the largest empirical gap from $\hat{\Delta}^k$-values 
is selected. On the other hand, if $1 \leq |A^k_f| \leq m_k$, the arm with the highest rank in $A_f^k$ in terms of $\hat{\mu}^k$-value, i.e., $\sigma_{\hat{\mu}^k}[A_f^k]^{-1}(1)$ is selected.
If $A_f^k$ is empty, then any arm in $A^k$ is selected.
After choosing $d_k$, CSAR sets $A^{k+1} = A^k \setminus \{d_k\}$.

Next, when $A_f^k$ is nonempty, CSAR checks the acceptance condition to determine whether $d_k$ also achieves the largest $\mu^k$-value among the arms.
If so, $d_k$ is accepted and $T_{m-m_k+1}$ is set to be $d_k$ and updates $m_k$ such that $m_{k+1} = m_k - 1$. If not, CSAR rejects $d_k$ (or does not set to a top-$m$ feasible arm).
On the other hand, when $A^k_f = \emptyset$, CSAR updates $m_k$ to count the ``empty" loop such that $m_{k+1} = m_k - 1$.
CSAR increases $k$ by 1 for the next phase and checks the exit condition. If $m_k=0$,
CSAR exits the loop and terminates with its output. Otherwise, it starts the next phase.

The pseudocode for the algorithm is given below.
\\
\\
\noindent\textbf{Constrained Successive Accept or Reject (CSAR) Algorithm}
\begin{itemize}
\item[1.] \textbf{Initialization:} Set $T_i=$Nil for $i=1,...,m$, $A^1 = A$, $m_1=m$, and $k=1$.\\
Select $N_i$ for $i=1,...,|A|$ such that $\sum_{k=1}^{|A|} N_k (|A|-k+1) \leq H$.
\item[2.] \textbf{Loop while $k\leq |A|$}
\begin{itemize}
\item[2.1] Play each $a \in A^k$ $N_k$ times and obtain $A^k_f = \{ a | \hat{C}^k(a) \leq \tau, a\in A^k\}$.
\item[2.2] \textbf{Deactivation:} $A^{k+1} = A^k \setminus \{d_k\}$, where
\begin{eqnarray*}
d_k = 
\begin{cases}
\sigma_{\hat{\Delta}^k}[A_f^k]^{-1}(1) & \mbox{if } |A^k_f| > m_k \\
\sigma_{\hat{\mu}^k}[A_f^k]^{-1}(1)  & \mbox{if } 1 \leq |A^k_f| \leq m_k \\
\mbox{any } a \in A^k & \mbox{if } |A^k_f| = 0.
\end{cases} 
\end{eqnarray*}
\item[2.3] \textbf{Accept or Reject:} If $A^k_f \neq \emptyset$ and 
$d_k = \sigma_{\hat{\mu}^k}[A_f^k]^{-1}(1)$, then $T_{m-m_k+1} = d_k$, $m_{k+1} = m_k - 1$. 
If $A^k_f = \emptyset$, $m_{k+1} = m_k - 1$.
\item[2.4] \textbf{Exit test:} Set $k \leftarrow k+1$. If $m_{k} = 0$, exit the loop.
\end{itemize}
\item[3.] \textbf{Output:} Set $T=\{T_i | T_i \neq \mbox{Nil}, i=1,...,m\}$. Return $T$.
\end{itemize}
\vspace{0.5cm}

We remark that one could consider a much simpler algorithm that at each phase $k$ with $A_f^k\neq \emptyset$, the empirically best arm $\sigma_{\hat{\mu}^k}[A_f^k]^{-1}(1)$ in $A_f^k$ is just selected as a top-$m$ feasible-arm.
(The arms selected at the previous phases have been all deactivated so that they are excluded from the current selection.)
This essentially corresponds to employing successively the method of sample average approximation (SAA)~\cite{kleywegt} over the feasible-arm set identified at the current phase. 
To analyze the probability of the correct identification, we would then consider the event that $A^k_f$ is correct and the selection by SAA is correct, i.e., $\sigma_{\hat{\mu}^k}[A_f^k]^{-1}(1)$ is in $\{\sigma_{\mu}[A^k_f]^{-1}(1),...,\sigma_{\mu}[A^k_f]^{-1}(m_k)\}$, where $m_k = m - (k-1)$ under the condition that all previous $k-1$ selections are correct.
Even if the probability can be bounded by an upper bound that converges to one but has the form of $(1-\exp(-H C))^m$ for a constant $C$ by $m$ independent correct-selections.
The \emph{exponential dependence on} $m$ makes the convergence slow (in particular as $m$ grows).
In CSAR, not only accepting correctly but also rejecting correctly at each phase is considered when the feasible active-arm set is nonempty.
In the step 2.3, CSAR accepts a largest empirical-gap achieving arm only if it is also the empirically best arm and rejects otherwise. This allows to derive a probability bound that does not depend exponentially on $m$ (see, Section~\ref{sec:pf}).

Another possible simpler approach than CSAR is to decompose the process of finding solution into the feasibility identification and the selection of the top feasible arms. For example, once we obtain $A_f^1$, we use (properly corrected) SAR to select the top arms in $A_f^1$.
This is plausible because if the correctness of $A_f^1$ is guaranteed in probability at some degree,
the selection process will provide the output with some guarantee in probability.
Analyzing the guarantee would necessarily consider an event that $A_f^1$ is correct and SAR makes correct identifications \emph{conditioned on the event that} $A_f^1$ is correct.
Because the overall convergence rate is affected by the probability that $A_f^1$ is correct, the convergence rate of this approach becomes slow.
Note that in this approach, any of the cost samples obtained for the phase $k\geq 2$ is \emph{not} utilized in the selection process.
Indeed, the similar idea was already employed in~\cite{changcmab} when designing a stochastic algorithm but the algorithm is shown to attain an exponential convergence rate only \emph{asymptotically}, i.e., for a sufficiently large value of $H$.
In CSAR, the feasibility identification and a top feasible-arm selection process are done together at the same phase \emph{by incorporating all of the samples available up to that phase}. This leads to an exponential convergence-rate of CSAR 
uniformly over (not just asymptotically) all finite horizons no smaller than $|A|$.

\section{Performance analysis}
\label{sec:pf}

Let $\mathcal{H}_k$ denote the $k$th harmonic number such that $\mathcal{H}_k = \sum_{i=1}^k 1/i, k\in Z^+$.
Set $N_k$ such that $N_k = n_k - n_{k-1}, k=1,2,...,|A|$, where $n_0=0$. Then trivially, $n_k$ is equal to the total number of plays for each active arm at the phase $k$. In other words, we have $n_k$ reward samples and $n_k$ cost samples for each active arm in $A^k$.
We consider a particular setting for our analysis:
\[ 
   n_k = \left \lceil \frac{H-|A|}{(|A|+1-k) \mathcal{H}_{|A|}} \right \rceil.
\] 
First, because $\sum_{k=1}^{|A|} (|A|-(k-1)) N_k = \sum_{k=1}^{|A|} n_k$, bounding the ceiling function in $n_k$ leads to $\sum_{k=1}^{|A|} n_k \leq \sum_{k=1}^{|A|} ( (H-|A|)/(|A|+1-k) \mathcal{H}_{|A|}^{-1} ) + 1$. The right hand side of the inequality is then less than equal to $(H-|A|) \mathcal{H}_{|A|} \mathcal{H}_{|A|}^{-1} + |A| = H.$ 
Therefore, the values of $N_k$'s make the total number of plays satisfy the budget $H$.
Second, $n_k$'s are monotonically decreasing. In other words, the more samples are allocated in the earlier phases. 
Finally, we will further see that this setting makes the contribution on the upper bound on the probability of the incorrect identification from the phase $k$ have the exponential form of $\exp(-n_k C)$ for a constant $C$.

We remark that a similar allocation methodology was already employed in SAR~\cite{bubeck} and a best-arm identification algorithm, ``successive rejects (SR)," in~\cite{audibert}. 
It is our purpose to maintain with proper changes the salient features of SAR, which extends SR.
In particular, Audibert \emph{et al.}~\cite{audibert} showed that their allocation rumake the total number of plays le yields an ``optimal" exponential convergence rate to zero for the probability of the incorrect identification by their algorithm. Even though we do not cover the optimality of the convergence rate in this paper, the method of setting $N_k$'s as above \emph{suffices our purpose}
of establishing an exponential convergence rate in the total number of plays limited by $H$ (see, e.g.,~\cite{shah} for other settings in allocation of sampling budget among arms and related algorithms).

The theorem below provides an upper bound on the probability of incorrect identification by CSAR. 
The bound approaches to zero exponentially fast in $H$ but the rate is compensated by the size of $A$ and a problem complexity constant, $\Delta_{\min}$, as given in the statement.
The assumption on the $\mu$-values and the $C$-values are necessary to have a meaningful result.
This can be relaxed by introducing some ``tolerance" parameters (cf., a remark in the concluding section).

\begin{thm}
\label{thm:main}
Assume that for all $a, a'\in A$, $\mu(a) \neq \mu(a')$ if $a\neq a'$ in $A$ and
$C(a) \neq \tau$ for all $a\in A$. Let $\Delta_{\min} = \min \{(\min_{a\in A} \Delta_c(a))^2/2, (\min_{a\in A} \Delta(a))^2/8\}$. 
Then the probability of the incorrect identification by CSAR is bounded above by
\[
2|A|^2 \exp \biggl (-\frac{(H-|A|)\Delta_{\min}}{|A|(\ln|A|+1)} \biggr ).
\]
\end{thm}
\vspace{0.5cm}

\begin{proof}
For each phase $k\in \{1,...,|A|\}$ at the beginning of the loop of CSAR,
let 
$\zeta_k$ denote the event
\[\biggl \{ A^k \neq \emptyset \wedge \forall a \in A^k \mbox{ } |\hat{C}^k(a) - C(a)| \leq \frac{1}{2}\Delta_c(a) \biggr \} 
\]
and 
$\chi_k$ denote the event 
\[\biggl\{ A_f^k \neq \emptyset \wedge \forall a \in A_f^k \mbox{ } |\hat{\mu}^k(a) - \mu(a)| \leq \frac{1}{4} \max_{a \in A_f^k} \Delta(a) \biggr\}.
\]

Because CSAR starts with $A^1=A$ and exactly one arm is deactivated, $A^k$ is nonempty
at each phase $k=1,...,|A|$ at the beginning of the loop.

\textbf{Identification of the set of feasible arms:}
Conditioned on $\zeta_k$, for any $a\in A^k$ such that $C(a) \leq \tau$,
$\hat{C}^k(a) - \tau \leq C(a) - \tau + \frac{1}{2} \Delta_c(a) = C(a) - \tau + \frac{1}{2}(\tau-C(a)) = \frac{1}{2} (C(a) - \tau) \leq 0.$
On the other hand, if $C(a) > \tau$ for $a\in A^k$, $\hat{C}^k(a) > \tau$ because
$\hat{C}^k(a) - \tau \geq C(a) - \tau - \frac{1}{2} \Delta_c(a) = C(a) - \tau - \frac{1}{2}(C(a) - \tau) = \frac{1}{2} (C(a) - \tau) > 0.$
This implies that $C(a) \leq \tau$ if and only if $\hat{C}^k(a) \leq \tau$. For any $k \in \{1,...,|A|\}$, on $\zeta_k$, $A_f^k$ is correct in that either it contains all feasible arms in $A^k$ or it is empty (all arms in $A^k$ are infeasible).
\\
\\
\noindent\textbf{Case $|A_f| \leq m \leq |A|$:}
Suppose that $A_f \neq \emptyset$. Condition on $\bigwedge_{k=1}^m\zeta_k$.
We first show that at $k=1,...,|A_f|$, $|A_f^k| = |A_f| - (k-1)$, $T_{k} =d_k$, and $m_{k+1} = m-k$ by induction on $k$.

For the base case, because $A_f^1$ is correct, $m_1=m \geq |A_f^1| = |A_f|$. CSAR
sets that $d_1 = \sigma_{\hat{\mu}^1}[A_f^1]^{-1}(1)$. $d_1$ is removed from $A^1$ and accepted
as $T_{m-m_k+1} = T_1$. We have that $T_1 = d_1$ and $m_2=m_1-1$.
Suppose that at $k$, the induction hypothesis is true. 
At $k+1$, because $A_f^{k+1}$ is correct, $|A_f^{k+1}| = |A_f^k| -1$ (since $d_k$ was removed from $A^k$).
Since $|A_f| \leq m$, we have that $m_{k+1} = m - k \geq |A_f| - (k-1) - 1$, $m_{k+1} \geq |A_f^{k+1}|$. By the choice of $d_{k+1}$ in the step 2.2, 
the acceptance condition is satisfied and CSAR sets $T_{m-m_{k+1} + 1} = T_{k+1} = d_{k+1}$ and $m_{k+2} = m_{k+1} -1 = (m-k)-1$.

The previous induction implies that at $k=|A_f|$, $T_i=d_i$ for $i=1,...,|A_f|$.
From $k=|A_f|+1,...,m$, because on $\zeta_k$, $A_f^k$ is correct, $A_f^k=\emptyset$. 
CSAR selects an arbitrary arm $a$ in $A^k$ as $d_k$ in the step 2.2, deactivates, and rejects this arm.
CSAR then sets $m_{k+1} = m_k -1$. By continuing this process, when $k=m$, $m_{k+1}=0$.
CSAR exits the loop and returns $T= \{T_1,...,T_{|A_f|}\}$ as the correct output.

On the other hand, if $A_f=\emptyset$, the event $\bigwedge_{k=1}^m\zeta_k$ implies that for $k=1,...,m$, $A_f^k=\emptyset$. Therefore CSAR sets $T=\emptyset$ as the correct output.

In sum, on $\bigwedge_{k=1}^m\zeta_k$, CSAR is correct when $|A_f| \leq m \leq |A|$.
\\
\\
\noindent\textbf{Case $1 \leq m < |A_f|$:}
For all $k=1,...,|A_f|-1$, 
conditioned on $\xi_k = \bigwedge_{i=1}^k (\{ m_k\geq 1 \} \wedge \zeta_k \wedge \chi_k)$, 
we show that CSAR either accepts correctly or reject correctly.

Fix $k \in \{2,...,|A_f|-1\}$. Conditioned on $\xi_{k}$, assume that for all $p=1,...,k-1$,
either CSAR accepted correctly such that $d_p$ was in $\{\sigma_{\mu}[A_f^p]^{-1}(1),...,\sigma_{\mu}[A_f^p]^{-1}(m_p)\}$
or CSAR rejected correctly such that $d_p$ was in $\{\sigma_{\mu}[A_f^p]^{-1}(m_p+1),...,\sigma_{\mu}[A_f^p]^{-1}(|A_f^p|)\}$.

For each $p=1,...,k-1$, let $q_p$ be the number of times CSAR accepted correctly up to $p$. Then $q_{k-1} \leq m-1$ so that $m_k \geq 1$ and $A^f_k$ is correct. Because $m < |A_f|$, $m_k = m-q_{k-1} < |A_f| - q_{k-1} = |A_f^k|$, therefore CSAR sets $d_k=\sigma_{\hat{\Delta}^k}[A^k_f]^{-1}(1)$ in the step 2.2.
\\
\\
\textbf{Correctness of acceptance:} 
Suppose that CSAR accepts $d_k$ but $d_k$ is not in $\{\sigma_{\mu}[A^k_f]^{-1}(1),...,\sigma_{\mu}[A^k_f]^{-1}(m_k)\}$.

Because $\hat{\mu}^k(\sigma_{\hat{\mu}^k}[A^k_f]^{-1}(1)) \geq \hat{\mu}^k(\sigma_{\mu}[A^k_f]^{-1}(1))$, $\chi_k$ implies that
\[
\mu(\sigma_{\hat{\mu}^k}[A^k_f]^{-1}(1)) + \frac{1}{4}\max_{a \in A_f^k} \Delta(a)  \geq \mu(\sigma_{\mu}[A^k_f]^{-1}(1)) - \frac{1}{4} \max_{a \in A_f^k} \Delta(a).
\]
Rearranging the terms lead to
\begin{eqnarray*}
\lefteqn{\max_{a \in A_f^k} \Delta(a) > \frac{1}{2} \max_{a \in A_f^k} \Delta(a)}\\
& & \geq 
\mu(\sigma_{\mu}[A^k_f]^{-1}(1)) - \mu(\sigma_{\hat{\mu}^k}[A^k_f]^{-1}(1)) \geq \mu(\sigma_{\mu}[A^k_f]^{-1}(1)) - \mu(\sigma_{\mu}[A^k_f]^{-1}(m+1)),
\end{eqnarray*} where the last inequality follows from the assumption that $\sigma_{\hat{\mu}^k}[A^k_f]^{-1}(1)$ is in $\{\sigma_{\mu}[A^k_f]^{-1}(m_k+1),...,\sigma_{\mu}[A^k_f]^{-1}(|A^k_f|)\}$. Therefore we have that
\begin{equation}
\label{1stcond}
  \max_{a \in A_f^k} \Delta(a) > \mu(\sigma_{\mu}[A^k_f]^{-1}(1)) - \mu(\sigma_{\mu}[A^k_f]^{-1}(m+1)).
\end{equation}

On the other hand, for all $i=1,...,m_k$,
\begin{eqnarray*}
\lefteqn{\hat{\mu}^k(\sigma_{\mu}[A^k_f]^{-1}(i)) \geq \mu(\sigma_{\mu}[A^k_f]^{-1}(i))  - \frac{1}{4} \max_{a \in A_f^k} \Delta(a)}\\
& & \geq \mu(\sigma_{\mu}[A^k_f]^{-1}(m_k))  - \frac{1}{4} \max_{a \in A_f^k} \Delta(a) \geq \mu(\sigma_{\mu}[A^k_f]^{-1}(m))  - \frac{1}{4} \max_{a \in A_f^k} \Delta(a),
\end{eqnarray*} where the last inequality holds because $m_k \leq m$. It follows that
\[
    \hat{\mu}^k(\sigma_{\hat{\mu}^k}[A^k_f]^{-1}(1)) \geq \hat{\mu}^k(\sigma_{\mu}[A^k_f]^{-1}(1)) \geq \mu(\sigma_{\mu}[A^k_f]^{-1}(m))  - \frac{1}{4} \max_{a \in A_f^k} \Delta(a).
\] Putting together, there exist $m_k+1$ arms $\{\sigma_{\mu}[A^k_f]^{-1}(1),...,\sigma_{\mu}[A^k_f]^{-1}(m_k), \sigma_{\hat{\mu}^k}[A^k_f]^{-1}(1)\}$ whose $\hat{\mu}^k$-values are bigger than equal to $\mu(\sigma_{\mu}[A^k_f]^{-1}(m))  - \frac{1}{4} \max_{a \in A_f^k} \Delta(a)$. This implies that
\begin{equation}
\label{mk}
    \hat{\mu}^k(\sigma_{\hat{\mu}^k}[A^k_f]^{-1}(m_k+1)) \geq \mu(\sigma_{\mu}[A^k_f]^{-1}(m))  - \frac{1}{4} \max_{a \in A_f^k} \Delta(a).
\end{equation}

By the acceptance condition of $\sigma_{\hat{\mu}^k}[A^k_f]^{-1}(1)  = \sigma_{\hat{\Delta}^k}[A^k_f]^{-1}(1)$,
\begin{equation} 
\label{accinq}
\hat{\Delta}^k(\sigma_{\hat{\mu}^k}[A^k_f]^{-1}(1)) \geq \hat{\Delta}^k(\sigma_{\hat{\mu}^k}[A^k_f]^{-1}(|A^k_f|)).
\end{equation}
We now upper and lower bounds the term in each side of~(\ref{accinq}), respectively:
\begin{eqnarray}
\lefteqn{\mu(\sigma_{\hat{\mu}^k}[A^k_f]^{-1}(1)) + \frac{1}{4} \max_{a \in A_f^k} \Delta(a) 
- \Bigl ( \mu(\sigma_{\mu}[A^k_f]^{-1}(m))  - \frac{1}{4} \max_{a \in A_f^k} \Delta(a) \Bigr )} \label{long1st}\\
& & \geq \hat{\mu}^k(\sigma_{\hat{\mu}^k}[A^k_f]^{-1}(1)) - \hat{\mu}^k(\sigma_{\hat{\mu}^k}[A^k_f]^{-1}(m_k+1)) \nonumber \mbox{ from~(\ref{mk}) and } \chi_k\\
& & = \hat{\Delta}^k(\sigma_{\hat{\mu}^k}[A^k_f]^{-1}(1)) \geq \hat{\Delta}^k(\sigma_{\hat{\mu}^k}[A^k_f]^{-1}(|A^k_f|)) \nonumber \mbox{ from~(\ref{accinq})}\\
& & = \hat{\mu}^k(\sigma_{\hat{\mu}^k}[A^k_f]^{-1}(m_k)) - \hat{\mu}^k(\sigma_{\hat{\mu}^k}[A^k_f]^{-1}(|A^k_f|)) \nonumber \mbox{ by the definition of } \hat{\Delta}^k\\
& & \geq \mu(\sigma_{\hat{\mu}^k}[A^k_f]^{-1}(m_k)) - \frac{1}{4} \max_{a \in A_f^k} \Delta(a) - \hat{\mu}^k(\sigma_{\hat{\mu}^k}[A^k_f]^{-1}(|A^k_f|)) \nonumber \mbox{ from } \chi_k\\
& & \geq \mu(\sigma_{\hat{\mu}^k}[A^k_f]^{-1}(m_k)) - \frac{1}{4} \max_{a \in A_f^k} \Delta(a) - \hat{\mu}^k(\sigma_{\mu}[A^k_f]^{-1}(|A^k_f|)) \nonumber \\
& & \geq \mu(\sigma_{\mu}[A^k_f]^{-1}(m)) - \frac{1}{4} \max_{a \in A_f^k} \Delta(a) - \Bigl (\mu(\sigma_{\mu}[A^k_f]^{-1}(|A^k_f|)) + \frac{1}{4} \max_{a \in A_f^k} \Delta(a) \Bigr ) \label{longlst}
\end{eqnarray} where the last inequality holds from $\chi_k$ and $m_k \leq m$.
Then rearranging the terms in~(\ref{long1st}) and~(\ref{longlst}), we have that
\begin{eqnarray}
\label{2ndcond}
\lefteqn{\max_{a \in A_f^k} \Delta(a)}\nonumber \\
& & \geq 2 \mu(\sigma_{\mu}[A^k_f]^{-1}(m)) - \mu(\sigma_{\hat{\mu}^k}[A^k_f]^{-1}(1)) - \mu(\sigma_{\mu}[A^k_f]^{-1}(|A^k_f|))  \nonumber \\
& & = \mu(\sigma_{\mu}[A^k_f]^{-1}(m)) - \mu(\sigma_{\mu}[A^k_f]^{-1}(|A^k_f|)) + \Bigl ( \mu(\sigma_{\mu}[A^k_f]^{-1}(m)) - \mu(\sigma_{\hat{\mu}^k}[A^k_f]^{-1}(1)) \Bigr ) \nonumber \\
& & > \mu(\sigma_{\mu}[A^k_f]^{-1}(m)) - \mu(\sigma_{\mu}[A^k_f]^{-1}(|A^k_f|)) 
\end{eqnarray} because $\mu(\sigma_{\mu}[A^k_f]^{-1}(m)) - \mu(\sigma_{\hat{\mu}^k}[A^k_f]^{-1}(1)) > 0$ from the assumption that $\sigma_{\hat{\mu}^k}[A^k_f]^{-1}(1)$ is not a top-$m$ arm.

Combining the results of~(\ref{1stcond}) and~(\ref{2ndcond}) leads to
\[
    \max_{a \in A_f^k} \Delta(a)  > \max\{\mu(\sigma_{\mu}[A^k_f]^{-1}(1)) - \mu(\sigma_{\mu}[A^k_f]^{-1}(m+1)),\mu(\sigma_{\mu}[A^k_f]^{-1}(m)) - \mu(\sigma_{\mu}[A^k_f]^{-1}(|A^k_f|)) \}.
\] This contradicts 
\[\max_{a \in A_f^k} \Delta(a) \leq \max\{\mu(\sigma_{\mu}[A^k_f]^{-1}(1)) - \mu(\sigma_{\mu}[A^k_f]^{-1}(m_k+1)),\mu(\sigma_{\mu}[A^k_f]^{-1}(m)) - \mu(\sigma_{\mu}[A^k_f]^{-1}(|A^k_f|)) \},
\] where this holds from the definition of $\Delta$ and the $\Delta$-values restricted to the elements of $A_f^k$. 
Therefore, CSAR accepts correctly such that $d_k$ is in $\{\sigma_{\mu}[A^k_f]^{-1}(1),...,\sigma_{\mu}[A^k_f]^{-1}(m_k)\}$.
\\
\\
\textbf{Correctness of rejection:} 
Because $d_k=\sigma_{\hat{\Delta}^k}[A^k_f]^{-1}(1)$ is rejected, observe that $\sigma_{\hat{\mu}^k}[A^k_f]^{-1}(|A^k_f|)  = d_k$.
Suppose that rejection is incorrect. Then $d_k$ is in $\{\sigma_{\mu}[A^k_f]^{-1}(1),...,\sigma_{\mu}[A^k_f]^{-1}(m_k)\}$. 

From $\hat{\mu}^k(\sigma_{\hat{\mu}^k}[A^k_f]^{-1}(|A_f^k|)) \leq \hat{\mu}^k(\sigma_{\mu}[A^k_f]^{-1}(|A_f^k|))$ and $\chi_k$, we have that 
\[
\max_{a \in A_f^k} \Delta(a) >
\mu(\sigma_{\hat{\mu}^k}[A^k_f]^{-1}(|A_f^k|)) - \mu(\sigma_{\mu}[A^k_f]^{-1}(|A_f^k|)).
\] Because $\mu(\sigma_{\hat{\mu}^k}[A^k_f]^{-1}(|A_f^k|)) \geq \mu(\sigma_{\mu}[A^k_f]^{-1}(m))$ from the assumption that $d_k$ is a top-$m$ arm,
\begin{equation}
\label{1stcondrej}
 \max_{a \in A_f^k} \Delta(a) > \mu(\sigma_{\mu}[A^k_f]^{-1}(m)) - \mu(\sigma_{\mu}[A^k_f]^{-1}(|A_f^k|)).
\end{equation}

On the other hand, all of the arms in $\{\sigma_{\mu}[A^k_f]^{-1}(m_k+1),...,\sigma_{\mu}[A^k_f]^{-1}(|A^k_f|)\}$ are non-top $m$ arms (i.e., a subset of $A_f \setminus \{\sigma_{\mu}[A^k_f]^{-1}(1),...,\sigma_{\mu}[A^k_f]^{-1}(m)\}$.
From $\chi_k$, for all $i=m_k+1,...,|A_f^k|$, 
\[
\hat{\mu}^k(\sigma_{\mu}[A^k_f]^{-1}(i)) \leq \mu(\sigma_{\mu}[A^k_f]^{-1}(i)) + \frac{1}{4} \max_{a \in A_f^k} \Delta(a)\leq \mu(\sigma_{\mu}[A^k_f]^{-1}(m+1))  + \frac{1}{4} \max_{a \in A_f^k} \Delta(a).
\] Furthermore,
\begin{equation}
\label{eqn:Ak}
    \hat{\mu}^k(\sigma_{\hat{\mu}^k}[A^k_f]^{-1}(|A_f^k|)) \leq \hat{\mu}^k(\sigma_{\mu}[A^k_f]^{-1}(m+1)) \leq \mu(\sigma_{\mu}[A^k_f]^{-1}(m+1))  + \frac{1}{4} \max_{a \in A_f^k} \Delta(a).
\end{equation} 
We see that there exist $|A^k_f| - m_k+1$ arms whose $\hat{\mu}^k$-values are less than equal to $\mu(\sigma_{\mu}[A^k_f]^{-1}(m+1))  + \frac{1}{4} \max_{a \in A_f^k} \Delta(a)$. Therefore,
\begin{equation}
\label{eqn:mk}
    \hat{\mu}^k(\sigma_{\hat{\mu}^k}[A^k_f]^{-1}(m_k)) \leq \mu(\sigma_{\mu}[A^k_f]^{-1}(m+1))  + \frac{1}{4} \max_{a \in A_f^k} \Delta(a).
\end{equation}

We now lower and upper bounds each side of the inequality $\hat{\Delta}^k(\sigma_{\hat{\mu}^k}[A^k_f]^{-1}(|A^k_f|)) > \hat{\Delta}^k(\sigma_{\hat{\mu}^k}[A^k_f]^{-1}(1))$ satisfied by the rejection condition with $\chi_k$:
\begin{eqnarray}
\lefteqn{\mu(\sigma_{\mu}[A^k_f]^{-1}(m+1)) + \frac{1}{4} \max_{a \in A_f^k} \Delta(a)
- \Bigl ( \mu(\sigma_{\hat{\mu}^k}[A^k_f]^{-1}(|A_f^k|))  - \frac{1}{4} \max_{a \in A_f^k} \Delta(a) \Bigr )} \label{rejlong1st}\\
& & \geq \hat{\mu}^k(\sigma_{\hat{\mu}^k}[A^k_f]^{-1}(m_k)) - \hat{\mu}^k(\sigma_{\hat{\mu}^k}[A^k_f]^{-1}(|A^k_f|)) \mbox{ from~(\ref{eqn:mk}) and } \chi_k \nonumber \\
& & = \hat{\Delta}^k(\sigma_{\hat{\mu}^k}[A^k_f]^{-1}(|A^k_f|)) > \hat{\Delta}^k(\sigma_{\hat{\mu}^k}[A^k_f]^{-1}(1)) \mbox{ by the rejection condition} \nonumber \\
& & = \hat{\mu}^k(\sigma_{\hat{\mu}^k}[A^k_f]^{-1}(1)) - \hat{\mu}^k(\sigma_{\hat{\mu}^k}[A^k_f]^{-1}(m_k+1)) \mbox{ by the definition of } \hat{\Delta}^k \nonumber \\
& & \geq \hat{\mu}^k(\sigma_{\mu}[A^k_f]^{-1}(1)) - \hat{\mu}^k(\sigma_{\hat{\mu}^k}[A^k_f]^{-1}(m_k+1)) \nonumber \\
& & \geq \mu(\sigma_{\mu}[A^k_f]^{-1}(1)) - \frac{1}{4} \max_{a \in A_f^k} \Delta(a) - \Bigl (\mu(\sigma_{\mu}[A^k_f]^{-1}(m+1)) + \frac{1}{4} \max_{a \in A_f^k} \Delta(a) \Bigr ), \label{rejlonglst}
\end{eqnarray} where the last inequality follows from $\chi_k$ and 
applying~(\ref{eqn:mk}) from $\hat{\mu}^k(\sigma_{\hat{\mu}^k}[A^k_f]^{-1}(m_k+1))
\leq \hat{\mu}^k(\sigma_{\hat{\mu}^k}[A^k_f]^{-1}(m_k))$.
Rearranging the terms of~(\ref{rejlong1st}) and~(\ref{rejlonglst}), we have that
\begin{eqnarray}
\label{2ndcondrej}
\lefteqn{\max_{a \in A_f^k} \Delta(a)}\nonumber \\
& & > \left( \mu(\sigma_{\mu}[A^k_f]^{-1}(1)) - \mu(\sigma_{\mu}[A^k_f]^{-1}(m+1)) \right ) + \left (\mu(\sigma_{\hat{\mu}^k}[A^k_f]^{-1}(|A^k_f|)) - \mu(\sigma_{\mu}[A^k_f]^{-1}(m+1)) \right ) \nonumber \\
& & > \mu(\sigma_{\mu}[A^k_f]^{-1}(1)) - \mu(\sigma_{\mu}[A^k_f]^{-1}(m+1)) \mbox{ because } d_k \mbox{ is assumed to be a top-}m \mbox{ arm}.
\end{eqnarray} 

Combining the results of (\ref{1stcondrej}) and~(\ref{2ndcondrej}) yields the inequality of
\[
    \max_{a \in A_f^k} \Delta(a)  > \max\{\mu(\sigma_{\mu}[A^k_f]^{-1}(1)) - \mu(\sigma_{\mu}[A^k_f]^{-1}(m+1)),\mu(\sigma_{\mu}[A^k_f]^{-1}(m)) - \mu(\sigma_{\mu}[A^k_f]^{-1}(|A^k_f|)) \},
\] which makes the same contradiction as before.

Combining the arguments of the correct acceptance and the correct rejection implies that CSAR 
either accepts correctly or rejects correctly when $p=k$.

We now consider the base case of $p=1$ conditioned on $\xi_1$. In this case, $A_f^1=A_f$ and $m_1=m \geq 1$. With the same reasoning as above, either CSAR accepts correctly or rejects correctly from $A_f^1$.

This concludes that for all $k=1,...,|A_f|-1$, on $\xi_{k}$, CSAR either accepts correctly or rejects correctly at $k$.
\\
\\
\noindent\textbf{Probability bound:} 
It remains to show that the upper bound of the probability of incorrect identification in the statement is correct
for each case.

\noindent\textbf{Case $0< |A_f| \leq m \leq |A|$:} Previously, we showed that CSAR is correct when $0< |A_f| \leq m \leq |A|$ on the event $\bigwedge_{k=1}^m\zeta_k$.
It follows that
\begin{eqnarray}
\lefteqn{\Pr (\{ T \neq \{T_1,....,T_{|A_f|}\} \}) \leq 1 - \Pr \biggl (\bigwedge_{k=1}^m \zeta_k \biggr )}  \\
& & \leq \sum_{k=1}^{m}\sum_{a\in A^k} \Pr \Bigl ( \Bigl \{ \Bigl |\hat{C}^k(a) - C(a) \Bigr | > \frac{1}{2}\Delta_c(a) \Bigr \} \Bigr ) \label{2nd} \\
& & \leq \sum_{k=1}^{|A|}\sum_{a\in A} 2\exp\Bigl (-2n_k (\Delta_c(a)/2)^2 \Bigr ) \label{3rd}\\
& & \leq 2|A|^2 \exp \Bigl (-\frac{(H-|A|) (\min_{a\in A} \Delta_c(a))^{2}}{2|A|(\ln|A|+1)} \Bigr ), \label{4th}
\end{eqnarray} where (\ref{2nd}) is by the union bound and (\ref{3rd}) by Hoeffding inequality~\cite{hoeff} and (\ref{4th}) follows because
\begin{eqnarray*}
\lefteqn{n_k \Bigl (\min_{a \in A} \Delta_c(a) \Bigr )^{2}} \\
& & \geq \frac{(H-|A|) (\min_{a\in A} \Delta_c(a))^{2}}{|A|\mathcal{H}_{|A|}} \geq \frac{(H-|A|) (\min_{a\in A} \Delta_c(a))^{2}}{|A|(\ln|A|+1)},
\end{eqnarray*} where we used $\mathcal{H}_{|A|} \leq \ln |A| + 1$ for $|A|\geq 1$.
\\
\\
\noindent\textbf{Case $|A_f| = 0$:} Similarly, this case has also shown to be correct on $\bigwedge_{k=1}^m\zeta_k$. Therefore,
\begin{eqnarray*}
\Pr (\{ T \neq \emptyset \}) \leq 2|A|^2 \exp \Bigl (-\frac{(H-|A|)(\min_{a\in A} \Delta_c(a))^{2}}{2(\ln|A|+\frac{1}{2})|A|} \Bigr ).
\end{eqnarray*}
\noindent\textbf{Case $1\leq m < |A_f|$:} Observe that $\xi_{|A_f|-1}$ implies that $m_{|A_f|-1} = 1$ and CSAR must accept correctly at $k=|A_f|-1$. 
Otherwise, it contradicts the previous statement of the correct acceptances or the correct rejections for $k=1,...,|A_f|-1$.
It follows that
\[
\Pr (\{ T \neq A_f(m) \} ) \leq 1-\Pr(\xi_{|A_f|-1}).
\]
We have that
\begin{eqnarray*}
\lefteqn{\Pr(\xi_{|A_f|-1}^c) \leq \sum_{k=1}^{|A_f|-1}\sum_{a\in A^k} \Pr \Bigl ( \Bigl \{ \Bigl |\hat{C}^k(a) - C(a) \Bigr | > \frac{1}{2}\Delta_c(a) \Bigr \} \Bigr )} \\
& & \hspace{4cm} + \sum_{k=1}^{|A_f|-1}\sum_{a\in A^k_f} \Pr \Bigl ( \Bigl \{ \Bigl |\hat{\mu}^k(a) - \mu(a) \Bigr | > \frac{1}{4} \max_{a \in A_f^k} \Delta(a) \Bigr \}\Bigr ) \\
& & \hspace{0.5cm}\leq \sum_{k=1}^{|A|}\sum_{a\in A} 2\exp\Bigl (-2n_k (\Delta_c(a)/2)^2 \Bigr ) + \sum_{k=1}^{|A|}\sum_{a\in A} 2\exp \Bigl (-2n_k (\max_{a\in A_f^k}\Delta(a)/4)^2 \Bigr ) \\
& & \hspace{0.5cm}\leq 2|A|^2 \exp \Bigl (-\frac{(H-|A|) (\min_{a\in A} \Delta_c(a))^{2}}{2|A|(\ln|A|+1)} \Bigr )
+ 2|A|^2 \exp \Bigl (-\frac{(H-|A|)(\min_{a\in A} \Delta(a))^{2}}{8|A|(\ln|A|+1)} \Bigr ).
\end{eqnarray*} 
\end{proof}

We remark that
even if CSAR has been designed with focusing on top-$m$ feasible-arm selections but not on ranking, it turns out that CSAR also produces asymptotically correct ranking such that $\mu(T_1) > \mu(T_2) > \cdots > \mu(T_m)$ when $1\leq m < |A_f|$.
We discuss an upper bound on the probability of the correct ranking that approaches to one as $H$ goes to infinity. However, the convergence speed might be slow because the bound essentially has an exponential dependence on $m$.

Let $P_{\mbox{\footnotesize{rank}}}$ be the probability of the event that $\mu(T_1) > \mu(T_2) > \cdots > \mu(T_m)$ for the output $\{T_1,...,T_m\}$ produced by CSAR.
Let $\rho$ be the event that CSAR accepted $T_{i}$ at $k_i$ for $i=1,...,m$. Then $P_{\mbox{\footnotesize{rank}}}$ is lower bounded as follows:
\[
P_{\mbox{\footnotesize{rank}}} \geq \prod_{i=1}^m 1-\Pr \Bigl ( \Bigl \{ 
\mu(\sigma_{\hat{\mu}^{k_i}}[A_f^{k_i}]^{-1}(1)) \neq \max_{a\in A_f^{k_i}} \mu(a) 
\Bigl | \rho \Bigr \} \Bigr ),
\] where
\begin{eqnarray*}
\lefteqn{\Pr \Bigl ( \Bigl \{ 
\mu(\sigma_{\hat{\mu}^{k_i}}[A_f^{k_i}]^{-1}(1)) \neq \max_{a\in A_f^{k_i}} \mu(a) 
\Bigl | \chi \Bigr \} \Bigr )}\\ 
&  & \leq \sum_{j=2}^{m_{k_i}} \Pr \Bigl (
\Bigl \{\hat{\mu}^{k_i}(\sigma_{\hat{\mu}^{k_i}}[A_f^{k_i}]^{-1}(1))  < \hat{\mu}^{k_i}(\sigma_{\hat{\mu}^{k_i}}[A_f^{k_i}]^{-1}(j)) \Bigr \} \Bigr ) \\
& & \leq \sum_{j=2}^{m_{k_i}} \Pr \Bigl (
\Bigl \{\hat{\mu}^{k_i}(\sigma_{\hat{\mu}^{k_i}}[A_f^{k_i}]^{-1}(1)) -  \hat{\mu}^{k_i}(\sigma_{\hat{\mu}^{k_i}}[A_f^{k_i}]^{-1}(j)) - (-\phi_{k_i}(j)) \geq \phi_{k_i}(j)) \Bigr \} \Bigr ) \\
& & \leq \sum_{j=2}^{m_{k_i}} \exp(- n_{k_i} \phi_{k_i}(j)^{2}),
\end{eqnarray*} where the exponential bound is obtained by Hoeffding inequality~\cite{hoeff} to the sum of $2n_{k_i}$ independent random variables taking values in $[0,1]$ or in $[-1,0]$ whose expectation is $-\phi_{k_i}(j) = - 
\mu(\sigma_{\hat{\mu}^{k_i}}[A_f^{k_i}]^{-1}(j)) +
\mu(\sigma_{\hat{\mu}^{k_i}}[A_f^{k_i}]^{-1}(1)).
$ It follows that $P_{\mbox{\footnotesize{rank}}}$ is lower bounded by
\[
    \prod_{i=1}^m \Bigl (1 - (m_{k_i} - 1) \exp(-n_{k_i} \min_{j=2,...,m_{k_i}} \phi_{k_i}(j) \Bigr ),
\] which approaches to one as $H$ goes to infinity.

\section{Concluding Remarks}

The performance bound in Theorem~\ref{thm:main} is given under the assumption that the reward means are all different among the feasible arms (as is common in the literature; see, e.g.,~\cite{bubeck}~\cite{zhang}~\cite{xiao}) and the difference between the cost mean and $\tau$ for each arm is not zero.
If this assumption does not hold, the bound becomes pointless because $\Delta_{\min}=0$. The problem complexity becomes infinite in that we would need infinite number of samplings to be able to differentiate the arms with the equal reward means or to identify the feasible arms that satisfy the equal constraint.
We can handle this case with introducing some ``tolerance" parameters $\epsilon >0$ and $\delta >0$ into $\Delta_c$ and $\Delta$, respectively. 
$\Delta_c$ is redefined such that $\Delta_c(a) = |C(a)-\tau| + \epsilon$ for $a\in A$. For $a\in A_f$, $\Delta(a) = \mu_{(\sigma_{\mu}[A_f](a))} - \mu_{(m+1)} +\delta$ if $\sigma_{\mu}[A_f](a) \leq m$, $\mu_{(m)} - \mu_{(\sigma_{\mu}[A_f](a))} + \delta$ otherwise. For all $a\in A\setminus A_f$, $\Delta(a) = \min_{a\in A_f} \Delta(a)$.
We then can say that $a \in A$ is $\epsilon$-feasible if $C(a)\leq \tau - \epsilon$ and $\epsilon$-infeasible if $C(a) > \tau + \epsilon$. It can be checked out that $|\hat{C}^k(a) - C(a)| \leq 1/2\Delta_c(a)$ implies that if $C(a)  \leq \tau - \epsilon$ then $\hat{C}^k(a) \leq \tau$ and if $C(a) > \tau + \epsilon$ then $\hat{C}^k(a) > \tau$.
Similarly, the zero gap between $\mu$-values can be offset with $\delta$.

In CSAR, the same allocation $N_k$ is given to each active arm in $A^k$ at each phase $k$ 
to estimate the reward mean and the cost mean. 
An adaptive sampling approach would be
interesting such that the sampling budget $N_k |A^k|$ given at the phase $k$ can
be differently used among the active arms in $A^k$, in which case the performance 
analysis is expected to become more complex.

\end{document}